\documentclass[12pt,german]{article}
\usepackage[hypertex,backref]{hyperref}
\usepackage{amscd,amsfonts,amssymb,amsthm,amscd,amsmath,latexsym}
\usepackage[T1]{fontenc}
\usepackage{german}
\usepackage{graphicx}
\usepackage[latin1]{inputenc}
\usepackage[all]{xy}
\usepackage[dvips]{color}
\usepackage{multicol}
\usepackage{pslatex} 
\usepackage[top=2cm, bottom=3cm,left=2.5cm,right=3.5cm]{geometry} %manuelle Seitenränder

\begin{titlepage}
\title{Cohomology of the Free Loop Space of a Complex Projective Space}
\author{Nora Seeliger}

\end{titlepage}

\begin{document}
\maketitle

Let $\Lambda (\mathbb{C}P^n)$ denote the free loop space of the complex projective space $\mathbb{C}P^n$, i. e. %$\mathbb{C}P^n$ is the projective space of the vector space $\mathbb{C}^{n+1}$ of dimension $n+1$ over the complex numbers $\mathbb{C}$ and $\Lambda(\mathbb{C}P^n)$ is 
the function space $\mathrm{map}(S^1,\mathbb{C}P^n)$ of unbased maps from a circle $S^1$ into $\mathbb{C}P^n$ topologized with the compact open topology. In this note we show that despite the fact that the natural fibration $\Omega(\mathbb{C}P^n)\hookrightarrow \Lambda(\mathbb{C}P^n)\stackrel{eval}{\longrightarrow}\mathbb{C}P^n$ has a cross section its Serre spectral sequence does not collapse: Here $eval$ is the evaluation map at a base point $*$ $\in \mathbb{C}P^n$. This result was originally proven in [1] by means of the Eilenberg-Moore spectral sequence as well as in [2] by using Sullivan's theory of minimal models. Moreover it can be concluded from [3] even though it is not explicitely stated there: These results were obtained by using the energy function as a Morse function and computing the homology by determing the infinte dimensional submanifolds of the free loop space given as the critical points of the Morse function. Here we only make use of an elementary computation with the Serre spectral sequence.\\

This work is part of the author's diploma thesis at the Georg-August-Universität Göttingen under the supervision of Prof. Dr. Thomas Schick. I would like to thank the Professors Fred Cohen, Paul Goerss, Ran Levi and Larry Smith for kind advice and enriching discussion on the topic.\\

Let us begin with the observation that $\Omega(\mathbb{C}P^n) \simeq S^1\times \Omega(S^{2n+1})$. To see this, consider the principal fibration $S^1\hookrightarrow S^{2n+1}\rightarrow \mathbb{C}P^n$ given by the Hopf map, with classifying map $\mathbb{C}P^n \hookrightarrow \mathbb{C}P^{\infty}$, and loop it. We obtain a principal fibration (up to homotopy) $\Omega(S^{2n+1})\hookrightarrow\Omega(\mathbb{C}P^n)\rightarrow S^1$. Recall that $\pi_1(\Omega(\mathbb{C}P^n))\cong\mathbb{Z}\cong\pi_1(\mathbb{C}P^{\infty})$ which implies the existence of a cross section. The result follows since a principal fibration with a cross section is trivial.\\

Since $I=[0,1]$ is contractible the inclusion of $\mathbb{C}P^n$ as the space of constant maps into $\mathrm{map}(I,\mathbb{C}P^n)$ is a deformation retract. Let $eval$: $\mathrm{map}(I,\mathbb{C}P^n)\rightarrow\mathbb{C}P^n\times\mathbb{C}P^n$ denote the evaluation of a map $w:I\rightarrow\mathbb{C}P^n$ at the two end points $0,1\in I$. This map is a fibration and the fibre over the base point $(*,*)\in\mathbb{C}P^n\times\mathbb{C}P^n$ is the usual space of based loops $\Omega(\mathbb{C}P^n)$ at the base point of $\mathbb{C}P^n$. The restriction of $eval$ to the constant loops, identified with $\mathbb{C}P^n$ is just the diagonal map $\Delta :\mathbb{C}P^n\rightarrow\mathbb{C}P^n\times\mathbb{C}P^n$. Putting these facts together yields the following commutative diagram of fibrations:\\
\[
\xymatrix@R=16pt@C=20pt{ {\Omega(\mathbb{C}P^n)}\ar[1,0]^{id}\ar[0,1]&{\Lambda(\mathbb{C}P^n)}\ar[1,0]^{exp}\ar[0,1]^{eval}&{\mathbb{C}P^n}\ar[1,0]^{\Delta}\\
{\Omega(\mathbb{C}P^n)}\ar[0,1]\restore&{Map(I,\mathbb{C}P^n)}\ar[0,1]^{eval}&{\mathbb{C}P^n\times\mathbb{C}P^n}\\
&&&&{}\save[]+<-6cm,0cm>*\txt<20pc>{\textit{Diagram 1}:  A Map of Fibrations}\\
} 
\]\\

in which $\overset{exp}{\rightarrow}$ is the map between function spaces induced by the exponential map $I\rightarrow S^1$.

Examine the Serre spectral sequence $\{E^r,d^r\}$ with coefficients in $\mathbb{Z}$ of the bottom fibration. One has $ E^2\cong H^*(\Omega(\mathbb{C}P^n))\otimes H^*(\mathbb{C}P^n\times\mathbb{C}P^n)\cong E[u]\otimes\Gamma[y]\otimes \frac{\mathbb{Z}[c_1]}{(c_1^{n+1})}\otimes\frac{\mathbb{Z}[c_2]}{(c_2^{n+1})}$
and \\$E^r\Rightarrow H^*(\mathbb{C}P^n;\mathbb{Z})$ where $E[-]$ is the exterior algebra, $\Gamma[-]$ denotes the divided power algebra, $|u|=1, |y|=2n, |c_1|=|c_2|=2$.  Let $v = c_1\otimes 1- 1\otimes c_2$ and $w=c_1\otimes 1$ and observe that they form a basis for $ H^2(\mathbb{C}P^n\times\mathbb{C}P^n)$. For the term $E^2$ one has the following picture:\\
\xymatrix@R=10pt@C=10pt{
&&&&&&&&&&&\\
2n&y\ar[-1,0]&&&&&&&&&&\\
&&&&&&&&&&&\\
1&u\ar@{-}[-2,0]\ar[1,2]^{d^2}&&{u\otimes v \atop u\otimes w}&&{\cdots}&&{\cdots}&&{\bullet}&&\\
\ar[0,11] &{\bullet}&&{v \atop w}&&&&&&{{w^n\atop w^n(w-v)} \atop {\vdots\atop (w-v)^n}}&&\\
0&\ar@{-}[-2,0]\restore&&2&&{\cdots}&&{\cdots}&&2n&&\\
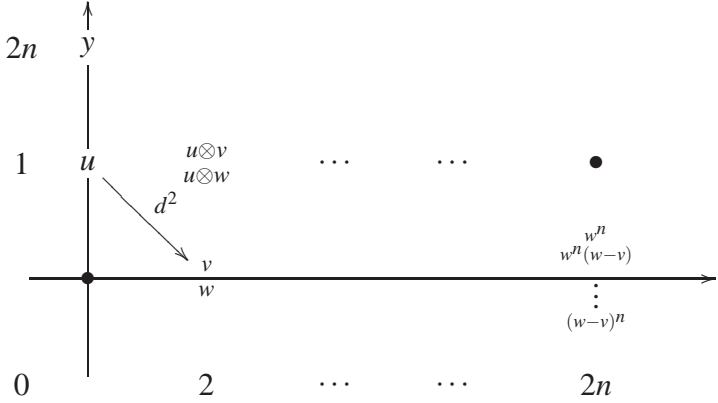
&&&&{}\save[]+<0cm,0cm>*\txt<20pc>{\textit{Diagram 2}: The Term $E^2$}\\
}\\[0.3cm]

The induced map $eval^*:H^*(\mathbb{C}P^n\times\mathbb{C}P^n)\rightarrow H^*(\mathrm{map}(I,\mathbb{C}P^n))$ coincides with the cup product $\Delta^*:H^*(\mathbb{C}P^n\times\mathbb{C}P^n)\rightarrow H^*(\mathbb{C}P^n)$.  Since  $eval^*$ is induced by the diagonal, we have $eval^*(v)=0$, so the element $v$ spans the kernel of this map. The only differential in the spectral sequence $\{ E^r, d^r \} $ whose target could be the element $v$
is the $d^2$-differential pictured in Diagram 2. This differential satisfies $d^2(y)\in E^2_{2,2n-1}$ and this homogeneous component of $E^2$ is zero unless $n=1$. So for $n>1$ we can compute the term $E^3$ as the homology of 
\begin{eqnarray*}
\{K,\delta\}=\{E[u]\otimes\frac{\mathbb{Z}[v,w]}{((v-w)^{n+1},w^{n+1})}, \delta(v)=0, \delta(u)=v\}\otimes\Gamma[y].
\end{eqnarray*}
The result is as pictured in the following diagram:\\
\xymatrix@R=10pt@C=10pt{
&&&&&&&&&\\
2n&y\ar[-1,0]\ar[2,6]^{d^{2n}}&&&&&&&&\\
&&&&&&&&&\\
1&&&&&&&{u\otimes (c_1^{n}\otimes 1+ c_1^{n-1}\otimes c_2+\cdots+1\otimes c_2^{n})}&&\\
\ar@{-}[0,3] &{\bullet}&&{w}\ar@{-}[0,4]&&&&{w^n}\ar[0,2]&&\\
0&\ar@{-}[-4,0]\restore&&2&&{\cdots}&&2n&&\\
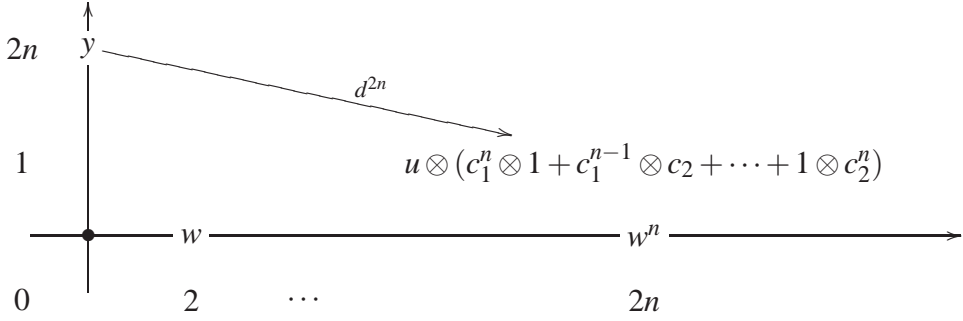
&&&&{}\save[]+<0cm,0cm>*\txt<20pc>{\textit{Diagram 3}: The Term $E^3$}\\
}\\[0.3cm]
The element $u\otimes (c_1^{n}\otimes 1+ c_1^{n-1}\otimes c_2+\cdots+1\otimes c_2^{n})$ has odd total degree, and the target of the spectral sequence $H^*(\mathbb{C}P^n)$ has no nonzero elements of odd degree. The only differential remaining in the spectral sequence that can annihilate $u\otimes (c_1^{n}\otimes 1+ c_1^{n-1}\otimes c_2+\cdots+1\otimes c_2^{n})$ is the $d^{2n}$ as pictured. It is also the first possible nonzero differential after $d^2$, so $E^3\cong E^{2n}$ and up to sign we must have $d^{2n}(y)=u\otimes (c_1^{n}\otimes 1+ c_1^{n-1}\otimes c_2+\cdots+1\otimes c_2^{n})$. The resulting $E^{2n+1}$ term then coincides with $H^*(\mathbb{C}P^n)$ so $E^{2n+1}=E^{\infty}$. This completes the description of the spectral sequence $\{E^r,d^r\}$. \\
Consider next the Serre spectral sequence $\{\overline{E}^r,\overline{d}^r\}$ for the free loop space fibration $\Omega(\mathbb{C}P^n)\hookrightarrow\Lambda(\mathbb{C}P^n)\rightarrow\mathbb{C}P^n$. The $\overline{E}^2$-term is as pictured in the following diagram:\\
\xymatrix@R=10pt@C=10pt{
&&&&&&&&&&&\\
2n&y\ar[-1,0]&&&&&&&&&&\\
&&&&&&&&&&&\\
1&u\ar[1,2]^{d^{2}}\ar@{-}[-2,0]&&{u\otimes x}&&&&{\cdots}&&{u\otimes x^n}&&\\
\ar@{-}[0,3] &{\bullet}&&x\ar@{-}[0,2]&&x^2\ar@{-}[0,4]&&&&{x^n}\ar[0,2]&&\\
0&\ar@{-}[-2,0]\restore&&2&&4&&{\cdots}&&2n&&\\
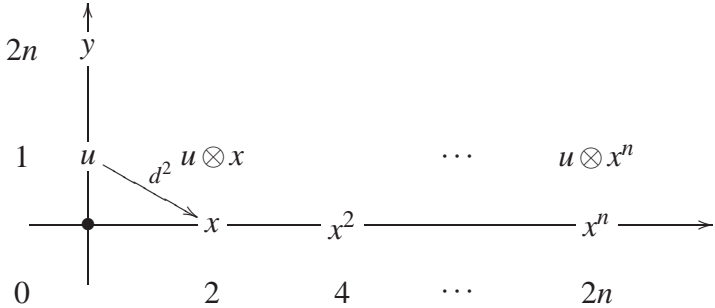
&&&&{}\save[]+<0cm,0cm>*\txt<20pc>{\textit{Diagram 4}: The Term $\overline{E}^2$ }\\
}\\[0.3cm]

This map of fibrations induces a map of Serre spectral sequences $\{E^r,d^r\}\stackrel{\phi}{\rightarrow}\{\overline{E}^r,\overline{d}^r\}$. Since $d^2(u)=v $ applying $\phi$ yields the formula
$\overline{d}^2(u)=\overline{d}^2(\phi (u))=\phi(d^2(u))=\phi (v)=0.$ So $\overline{d}^2(y)=0$ and hence for placement reasons one obtains $\overline{E}^2\cong\cdots\cong\overline{E}^{2n}$. The resulting term is pictured in the following diagram:\\[0.3 cm]
\xymatrix@R=10pt@C=10pt{
&&&&&&&&&&&\\
2n&y\ar[-1,0]\ar[2,8]^{d^{2n}}&&&&&&&&&&\\
&&&&&&&&&&&\\
1&u\ar@{-}[-2,0]&&{u\otimes x}&&&&{\cdots}&&{u\otimes x^n}&&\\
\ar@{-}[0,3] &{\bullet}&&x\ar@{-}[0,6]&&&&&&{x^n}\ar[0,2]&&\\
0&\ar@{-}[-2,0]\restore&&2&&&&{\cdots}&&2n&&\\
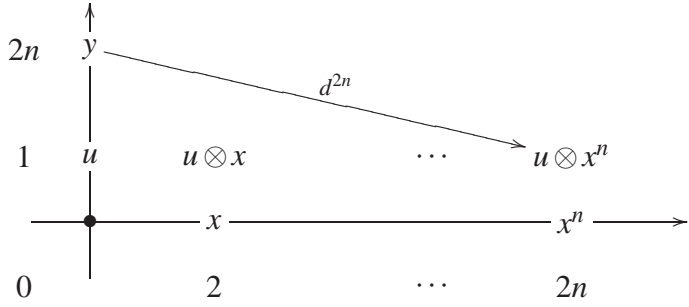
&&&&{}\save[]+<0cm,0cm>*\txt<20pc>{\textit{Diagram 5}:     The Term $\overline{E}^{2n}$}\\
}\\[0.3cm]

Again we can apply our knowledge of the spectral sequence $\{E^r,d^r\}$ and the map $\phi$ of spectral sequences to obtain:
\begin{eqnarray*}
\overline{d}^{2n}(y)=\overline{d}^{2n}(\phi(y))=\phi(d^{2n}(y))=\phi(u\otimes (c_1^{n}\otimes 1+ c_1^{n-1}\otimes c_2+\cdots+1\otimes c_2^{n}))=u\otimes (n+1)x^n.
\end{eqnarray*}
This element being nonzero we have proven the following result:\\[0.3 cm]
\textbf{Theorem:} The free loop space fibration $\Omega \mathbb{C}P^n\hookrightarrow\Lambda(\mathbb{C}P^n)\rightarrow \mathbb{C}P^n$ has a cross section and for $n \geq 1$ its integral Serre spectral sequence does not collapse.\\[0.3 cm]
Remarks: The theorem holds for n=1 since $\mathbb{C}P^1\cong S^2$. This can easily be seen by inspecting the corresponding Serre spectral sequences. The argument also goes through for coefficients in a field of characteristic $p$ as long as $p\nmid  n+1$. Moreover, the theorem holds for every space $X$ with $H^*(X;\mathbb{Q})\cong\frac{\mathbb{Q}[x]}{(x^{k})}$ and $H^*(\Omega(X);\mathbb{Q})\cong E[u]\otimes \mathbb{Q}[y]$ where $|x|=2m,|u|=2m-1$ and $|y|=2mk$. This includes $\mathbb{H}P^n$ as well as $S^{2k}$ and the Cayley projective plane. Rather little is known about the relation between $H^*(X)$ and $H^*(\Omega(X))$. Therefore it is difficult to formulate a generalization of our results. \\[0.3 cm]
This theorem appears under a totally different point of view as Lemma 2.14 in [5]. They don't generalise to the other globally symmetric spaces of rank 1. The cohomology of the free loop space of a complex projective space is discussed as well in [4].\\[0.3 cm]
\textbf{References:}\\
$[1]$ L. Smith: \textit{On the characteristic zero cohomology of the free loop space}, Amer. J. Math., vol. 103 (1981), pp. 887-910.\\
$[2]$ M. Vigué-Poirrier: \textit{Dans le fibré de l'espace des lacets libres, la fibre n'est pas, en général, totalement non cohomologue à zéro}, Math. Z. 181 (1982), pp. 537-542\\
$[3]$ W. Ziller: \textit{The Free Loop Space of Globally Symmetric Spaces} Inv. Math. 41, 1-22 (1977)\\
$[4]$ M. Crabb and I. James: \textit{Fiberwise homotopy theory}, Springer Monographs in Mathematics, Springer-Verlag London Limited, 1998. \\
$[5]$ E. Fadell, S. Husseini: \textit{Infinite cup length in free loop spaces with an application to a problem of the N-body type}, Ann. Institut Henri Poincare, Analyse non lineaire, vol. 9, n.3 (1992), p. 305 - 319.\\[0.3 cm]
Department of Mathematical Sciences, University of Aberdeen, Meston Building 401,\\ AB24 3UE Aberdeen, U. K,
e-mail: \textit{n.seeliger@maths.abdn.ac.uk}. 
\end{document}